\documentclass[a4paper, 12pt]{article}
 
\usepackage[top=2cm, bottom=2cm, left=2cm, right=2cm]{geometry}

\usepackage{graphicx, array, amsmath, amsthm, mathrsfs, framed, stmaryrd, marvosym, bbm, amssymb, mathtools, verbatim, tikz, tikz-cd}
\usepackage[all]{xy}

\tikzset{node distance=2cm, auto}

\usepackage{enumerate}

\usepackage{color}

\newtheorem{proposition}{Proposition}[section]
\newtheorem{lemma}[proposition]{Lemma}
\newtheorem{theorem}[proposition]{Theorem}

\newtheorem{notation}[proposition]{Notation}
\newtheorem{remark}[proposition]{Remark}

\theoremstyle{definition}
\newtheorem{definition}[proposition]{Definition}

\newcommand{\id}{\text{id}}

\newcommand{\tcr}{\textcolor{red}}
\newcommand{\co}{\colon}

\newcommand{\bs}{\backslash}
\newcommand{\be}{\begin{equation}}
\newcommand{\ee}{\end{equation}}
\newcommand{\piecewise}[5]{ #1 = \begin{cases} #2 & #3 \\ #4 & #5 \end{cases} }

\newcommand{\K}{\mathbb{K}}

\newcommand{\mc}[1]{\mathcal{#1}}

\newcommand{\ol}{\overline}
\newcommand{\ul}{\underline}

\newcommand{\g}{\mathfrak{g}}
\newcommand{\h}{\mathfrak{h}}
\newcommand{\n}{\mathfrak{n}}

\newcommand{\tx}{\text}

\newcommand{\ts}{\textstyle}

\newcommand{\N}{\mathbb{N}}
\newcommand{\C}{\mathbb{C}}

\newcommand{\Z}{\mathbb{Z}}
\newcommand{\Q}{\mathbb{Q}}

\newcommand{\BD}{\begin{definition}}
\newcommand{\BT}{\begin{theorem}}
\newcommand{\BL}{\begin{lemma}}
\newcommand{\ED}{\end{definition}}
\newcommand{\ET}{\end{theorem}}
\newcommand{\EL}{\end{lemma}}
\newcommand{\BProp}{\begin{proposition}}
\newcommand{\EProp}{\end{proposition}}
\newcommand{\BRem}{\begin{remark}}
\newcommand{\ERem}{\end{remark}}

\newcommand{\ang}[1]{\langle #1 \rangle}

\newcommand{\llb}{\llbracket}
\newcommand{\rrb}{\rrbracket}

\newcommand{\mf}{\mathfrak}

 {\begin{list}{}%
         {\setlength{\leftmargin}{13mm}}%
         \item[]%
 }
 {\end{list}}

 {\begin{list}{}%
         {\setlength{\leftmargin}{13mm}}%
         \item[]%
 }
 {\end{list}}

\DeclareMathOperator*{\rk}{rk}

\DeclareMathOperator*{\Prim}{Prim}
\DeclareMathOperator*{\Spec}{Spec}

\DeclareMathOperator*{\Frac}{Frac}

\DeclareMathOperator*{\trdeg}{tr.deg}

\DeclareMathOperator*{\height}{ht}

\DeclareMathOperator*{\locdeg}{loc.deg}
\DeclareMathOperator*{\primdeg}{prim.deg}
\DeclareMathOperator*{\ratdeg}{rat.deg}
\DeclareMathOperator*{\kdim}{K.dim}

\DeclareMathOperator*{\GKdim}{GKdim}
\DeclareMathOperator*{\Ad}{Ad}

\title{A strong Dixmier-Moeglin equivalence for quantum Schubert cells}
\author{Jason Bell, St\'ephane Launois, Brendan Nolan}
\date{\today}
\begin{document}
 
\maketitle

\begin{abstract}
\noindent Dixmier and Moeglin gave an algebraic condition and a topological condition for recognising the primitive ideals among the prime ideals of the universal enveloping algebra of a finite-dimensional complex Lie algebra; they showed that the primitive, rational, and locally closed ideals coincide. In modern terminology, they showed that the universal enveloping algebra of a finite-dimensional complex Lie algebra satisfies the \emph{Dixmier-Moeglin equivalence}. 

We define quantities which measure how ``close" an arbitrary prime ideal of a noetherian algebra is to being primitive, rational, and locally closed; if every prime ideal is equally ``close" to satisfying each of these three properties, then we say that the algebra satisfies the \emph{strong Dixmier-Moeglin equivalence}. Using the example of the universal enveloping algebra
of $\mf{sl}_2(\C)$, we show that the strong Dixmier-Moeglin equivalence is strictly stronger than the Dixmier-Moeglin equivalence.

For a simple complex Lie algebra $\g$, a non root of unity $q\neq 0$ in an infinite field $\K$, and an element $w$ of the Weyl group of $\g$, De Concini, Kac, and Procesi have constructed a subalgebra $U_q[w]$ of the quantised enveloping $\K$-algebra $U_q(\g)$. 
These \emph{quantum Schubert cells} are known to satisfy the Dixmier-Moeglin equivalence and we show that they in fact satisfy the strong Dixmier-Moeglin equivalence. Along the way, we show that commutative affine domains, uniparameter quantum tori, and uniparameter quantum affine spaces satisfy the strong Dixmier-Moeglin equivalence. 
\end{abstract}

\section{Introduction}
Throughout this paper, $\K$ denotes an infinite field of arbitrary characteristic and, unless otherwise stated, every algebra is a unital associative $\K$-algebra and every ideal is two-sided.

It is a difficult and often intractable problem to classify the irreducible representations of an algebra. Dixmier proposed that a good first step towards tackling this problem would be to find the kernels of the irreducible representations, that is the annihilators of the simple modules, namely the primitive ideals. 
 In any ring, every primitive ideal is prime; Dixmier \cite{Dixmier} and Moeglin \cite{Moeglin} gave an algebraic condition and a topological condition for deciding whether or not a given prime ideal of the universal enveloping algebra of a finite-dimensional complex Lie algebra is primitive:
\begin{itemize}
\item A prime ideal $P$ of a ring $R$ is said to be \emph{locally closed} if the singleton set $\{P\}$ is locally closed in the Zariski topology on $\Spec R$. Equivalently, $\{P\}$ is the intersection of a Zariski-open subset of $\Spec R$ and a Zariski-closed subset of $\Spec R$. (For a prime ideal $P$ of a ring $R$, it is easily shown that $P$ is locally closed if and only if $P$ is strictly contained in the intersection of all prime ideals of $R$ which strictly contain $P$.)
\item A prime ideal $P$ of a noetherian $\K$-algebra $R$ is said to be \emph{rational} if the field extension $\mc{Z}(\Frac R/P)$ of $\K$ is algebraic. 
\end{itemize}
Dixmier and Moeglin proved that for a prime ideal of the universal enveloping algebra of a finite-dimensional complex Lie algebra, the properties of being primitive, locally closed, and rational are equivalent. In modern terminology, they proved that the universal enveloping algebra of a finite-dimensional complex Lie algebra satisfies the \emph{Dixmier-Moeglin equivalence}.

Since the work of Dixmier and Moeglin on universal enveloping algebras of finite-dimensional complex Lie algebras, many more algebras have been shown to satisfy the Dixmier-Moeglin equivalence: \cite[Corollary II.8.5]{BG} lists several quantised coordinate rings which satisfy the Dixmier-Moeglin equivalence; the first named author, Rogalski, and Sierra \cite{BRS} have shown that twisted homogeneous coordinate rings of projective surfaces satisfy the Dixmier-Moeglin equivalence. However, Irving \cite{Irving} and Lorenz \cite{Lorenz} have shown that there exist noetherian algebras for which the Dixmier-Moeglin equivalence fails.

Our goal is to extend the notion of the Dixmier-Moeglin equivalence to all prime ideals, in a way which captures how ``close'' they are to being primitive.  Of course, not all non-primitive prime ideals are created equal.  For example, in the polynomial ring $\C[x,y]$, the primitive ideals are the maximal ideals $\ang{x-\alpha,y-\beta}$.  For this reason, we think of the prime ideal $\ang{x}$ as being ``closer" to being primitive than the prime ideal $\ang{0}$, in the same sense that it is ``closer" to being maximal --- that is, the height of $\ang{x}$ is greater than the height of $\ang{0}$.  

In general, given a noetherian $\K$-algebra $R$ and given a prime ideal $P$ of $R$, we are interested in the \emph{primitivity degree}, $\primdeg P$, of $P$, which we define as follows:
\[
\primdeg P:=\inf \{ \height Q\ | \ Q\in \Prim R/P \},
\] 
where $\Prim R/P$ denotes the subspace of $\Spec R/P$ consisting of the primitive ideals of $R/P$. 
This quantity gives a measure of how close the prime ideal $P$ is to being primitive. Clearly, $P$ is primitive if and only if $\primdeg P=0$.

We use this idea to extend the notion of the Dixmier-Moeglin equivalence to all prime ideals.  To this end, we define generalisations of the notions of a locally closed ideal and a rational ideal.

It is easy to extend the notion of a rational ideal: for a prime ideal $P$ of $R$, we define the \emph{rationality degree}, $\ratdeg P$, of $P$ to be the transcendence degree of the field extension $\mc{Z}(\Frac R/P)$ of $\K$. Clearly, $P$ is rational if and only if $\ratdeg P=0$. 

In the same spirit of generalisation, we define the \emph{local closure degree}, $\locdeg P$, of a prime ideal $P$ of $R$ to be the smallest nonnegative integer $d$ such that $\bigcap_{Q\in \Spec_{>d}R/P}Q\neq 0$, where $\Spec_{>d}R/P$ denotes the subspace of $\Spec R/P$ consisting of all prime ideals of $R/P$ which are of height strictly greater than $d$. Clearly, $P$ is locally closed if and only if $\locdeg P=0$. 
\begin{remark}
In the case that the noetherian $\K$-algebra $R$ has finite Gelfand-Kirillov dimension, all prime ideals of $R$ have finite height by \cite[Corollary 3.16]{KrauseLenagan}. 
All of the algebras which will concern us in this paper have finite Gelfand-Kirillov dimension and so we shall always use the following equivalent characterisation of local closure degree: for a prime ideal $P$ of $R$, $\locdeg P$ is the smallest nonnegative integer $d$ such that
$\bigcap_{Q\in \Spec_{d+1}R/P}Q\neq 0$, where $\Spec_{d+1}R/P$ denotes the subspace of $\Spec R/P$ consisting of all prime ideals of $R/P$ which are of height $d+1$. 
\end{remark}

\begin{definition}
A noetherian $\K$-algebra $R$ is said to satisfy the \emph{strong Dixmier-Moeglin equivalence} if every prime ideal $P$ of $R$ satisfies $\locdeg P=\primdeg P=\ratdeg P$.
\end{definition}

We remark that the strong Dixmier-Moeglin equivalence is strictly stronger than the Dixmier-Moeglin equivalence. Indeed the Dixmier-Moeglin equivalence simply says that if $P$ is a prime ideal of a noetherian $\K$-algebra $R$, then
\[
\locdeg P=0\iff \primdeg P=0\iff \ratdeg P=0.
\]
Even though the universal enveloping algebra, $U(\mathfrak{sl}_2(\C))$, of $\mathfrak{sl}_2(\C)$ satisfies the Dixmier-Moeglin equivalence (as was shown in the original work of Dixmier and Moeglin), it fails to satisfy the strong Dixmier-Moeglin equivalence. Indeed, since $U(\mathfrak{sl}_2(\C))$ is a domain, $\ang{0}$ is a (completely) prime ideal of  $U(\mathfrak{sl}_2(\C))$. By \cite[Remark 4.6]{Catoiu}, all nonzero prime ideals of $U(\mathfrak{sl}_2(\C))$ are primitive, so that $\primdeg\ang{0}=1$. It is well known that the centre of $U(\mathfrak{sl}_2(\C))$ is given by the polynomials in the Casimir element; by \cite[Corollary 4.2.3]{DixmiersBook}, $\mc{Z}(\Frac U(\mathfrak{sl}_2(\C)))$ is given by the rational functions in the Casimir element, so that $\ratdeg \ang{0}=\trdeg_\C \mc{Z}(\Frac U(\mathfrak{sl}_2(\C)))=1$. By \cite[Theorem 4.5 and Proposition 5.13]{Catoiu}, there are infinitely many height two prime ideals in $U(\mathfrak{sl}_2(\C))$ and their intersection is zero, so that $\locdeg \ang{0}>1$. Since, by \cite[Theorem 4.5]{Catoiu}, there are no height three prime ideals in $U(\mathfrak{sl}_2(\C))$, the intersection of the height three prime ideals is nonzero (in fact it is the entirety of $U(\mathfrak{sl}_2(\C))$), so that $\locdeg\ang{0}=2$. 

\vspace{2mm}

The goal of this paper is to prove that quantum Schubert cells, which we now briefly discuss (see Section \ref{section:q-schubert} for more details), satisfy the strong Dixmier-Moeglin equivalence. Let $\g$ be a simple complex Lie algebra of rank $n$ and let $\pi=\{\alpha_1,\ldots,\alpha_n\}$ be the set of simple roots associated to a triangular decomposition $\g=\n^-\oplus \h\oplus \n^+$ of $\g$.
 Where $q\in \K^\times$ is not a root of unity and $w$ is an element of the Weyl group of $\g$, De Concini, Kac, and Procesi \cite{DKP} 
defined a quantum analogue, $U_q[w]$, of the universal enveloping algebra of the nilpotent Lie algebra $\n^+\cap \Ad_w(\n^-)$. These \emph{quantum Schubert cells} $U_q[w]$ shall be our main objects of study.

\vspace{2mm}

It shall be useful to define a weaker version of the strong Dixmier-Moeglin equivalence which is often easy to prove and provides a useful stepping-stone to proving the strong Dixmier-Moeglin equivalence.
\begin{definition}
A noetherian $\K$-algebra $R$ is said to satisfy the \emph{quasi strong Dixmier-Moeglin equivalence} if every prime ideal $P$ of $R$ satisfies $\locdeg P=\ratdeg P$. 
\end{definition}

With the quasi strong Dixmier-Moeglin equivalence in hand for a noetherian $\K$-algebra $R$, the problem is reduced to showing that every prime ideal $P$ of $R$ satisfies $\primdeg P=\ratdeg P$. For a quantum Schubert cell $U_q[w]$, we prove this by exploiting the good behaviour of the poset of $H$-invariant prime ideals of $U_q[w]$, where $H$ is a suitable algebraic $\K$-torus acting rationally on $U_q[w]$ by $\K$-algebra automorphisms.

\vspace{2mm}

This paper is organised as follows. First, we prove various general results about the (quasi) strong Dixmier-Moeglin equivalence (Section \ref{section:generalities}). Next, we consider various examples from the quantum world. Using Cauchon's theory of deleting derivations, one can relate the prime and primitive spectra of a quantum Schubert cell to those of an associated uniparameter quantum affine space, which can in turn be related via localisations to the prime and primitive spectra of a family of uniparameter quantum tori. Since there is a bi-increasing homeomorphism between the prime spectrum of a uniparameter quantum torus and the prime spectrum of its centre, which is a commutative affine domain, we are guided into a natural strategy: we shall prove the strong Dixmier-Moeglin equivalence first for commutative affine domains (Section \ref{section:commutative}), then for uniparameter quantum tori (Section \ref{section:q-tori}), then for uniparameter quantum affine spaces (Section \ref{section:q-affine}), and finally for quantum Schubert cells (Section \ref{section:q-schubert}). Partial results are also obtained for a larger class of algebras --- we prove in Section \ref{CGL DDA} that every uniparameter Cauchon-Goodearl-Letzter extension satisfies the quasi strong Dixmier-Moeglin equivalence.

We have partial results for quantised coordinate rings and quantum Grassmannians and we have reason to believe that they satisfy the strong Dixmier-Moeglin equivalence; we will return to these algebras in a later paper.

\section{General results on the (quasi) strong Dixmier-Moeglin equivalence}
\label{section:generalities}

In this section we prove that, under some mild assumptions, the primitivity degree of a prime ideal is bounded above by its local closure degree, and then we prove transfer results for the quasi strong Dixmier-Moeglin equivalence for an algebra and its localisations.

\subsection{An upper bound for the primitivity degree}

Some of the implications needed to prove the Dixmier-Moeglin equivalence hold in a very general setting. Recall that a noetherian $\K$-algebra $R$ is said to satisfy the \emph{noncommutative Nullstellensatz} over $\K$ if $R$ is a Jacobson ring and the endomorphism ring of every irreducible $R$-module is algebraic over $\K$. By \cite[Lemma II.7.15]{BG}, for any noetherian $\K$-algebra $R$ which satisfies the noncommutative Nullstellensatz over $\K$ and for any prime ideal $P$ of $R$, we have 
\begin{equation}\label{eq implications}
P\tx{ is locally closed }\implies P\tx{ is primitive }\implies P\tx{ is rational.}
\end{equation}
We have generalised the first implication above to a large class of algebras:
\begin{proposition}\label{general ineq}
Let $R$ be a catenary, noetherian $\K$-algebra of finite Gelfand-Kirillov dimension which satisfies the noncommutative Nullstellensatz over $\K$. Then for any prime ideal $P$ of $R$, we have $\locdeg P\geq \primdeg P$.
\begin{proof}
We prove the result by induction on $d:= \locdeg P$. If $d=0$, then the result follows immediately from \cite[Lemma II.7.15]{BG}.

Suppose that $d>0$ and set $\ol{R}=R/P$. $\Spec_{d}\ol{R}$ is a set of nonzero prime ideals of $\ol{R}$ whose intersection is zero. It follows that the set of height one prime ideals of $\ol{R}$ has zero intersection. Hence there is a height one prime ideal $T$ of $\ol{R}$ which fails to contain the nonzero ideal $I:=\bigcap_{Q\in \Spec_{d+1}\ol{R}}Q$ of $\ol{R}$. 

Each prime ideal of $\ol{R}/T$ whose height is $d$ corresponds to a prime ideal of $\ol{R}$ which contains $T$, has height $d+1$, and hence contains $I+T$. 

Now each prime ideal of $\ol{R}/T$ whose height is $d$ must contain the image of $I+T$ in $\ol{R}/T$, which is nonzero. Hence $\bigcap_{Q\in \Spec_d \ol{R}/T}Q\neq 0$, so that $\locdeg T\leq d-1$. It follows from the induction hypothesis that $\primdeg T\leq d-1$. So there is a primitive ideal of $\ol{R}/T$ whose height is at most $d-1$; this corresponds to a primitive ideal of $\ol{R}=R/P$ whose height is at most $d$. Now $\primdeg P\leq d=\locdeg P$, as required. 
\end{proof}
\end{proposition}
We do not know whether the second implication in \eqref{eq implications} can be similarly generalised but we will prove, on a case-by-case basis, that for a prime ideal $P$ of a commutative affine domain, a uniparameter quantum torus, a uniparameter quantum affine space, or a quantum Schubert cell, we have 
\[
\primdeg P=\ratdeg P.
\]

\subsection{Transferring the quasi strong Dixmier-Moeglin equivalence}
Recall that a noetherian $\K$-algebra $R$ is said to satisfy the quasi strong Dixmier-Moeglin equivalence if, for every prime ideal $P$ of $R$, we have $\locdeg P=\ratdeg P$. 
\begin{lemma}\label{localise and loc deg} 
Let $R$ be a noetherian $\K$-algebra of finite Gelfand-Kirillov dimension which is a domain and in which every prime ideal is completely prime. Let $\mc{E}$ be a right Ore set of regular elements of $R$ which is finitely generated as a multiplicative system. Then for any $d\in \N\bs\{0\}$, we have 
\[
\ts{ \bigcap_{P\in \Spec_d R}P\neq 0\iff \bigcap_{Q\in \Spec_d R\mc{E}^{-1}}Q\neq 0. }
\]
It follows immediately that $\locdeg \ang{0}_R=\locdeg \ang{0}_{R\mc{E}^{-1}}$, where $\ang{0}_R$ and $\ang{0}_{R\mc{E}^{-1}}$ denote the zero ideals of $R$ and $R\mc{E}^{-1}$ respectively.  
\begin{proof}
Let $\mc{E}$ be generated as a multiplicative system by $x_1,\ldots,x_n$. Since all prime ideals of $R$ are completely prime, the conditions $P\cap \mc{E}=\emptyset$ and $x_1,\ldots,x_n\notin P$ are equivalent for every prime ideal $P$ of $R$.

By \cite[Theorem 10.20]{GW}, extension ($P\mapsto P\mc{E}^{-1}$) and contraction ($Q\mapsto Q\cap R$) are mutually inverse increasing homeomorphisms between $\{P\in \Spec R \ | \ P\cap \mc{E}=\emptyset\}=\{P\in \Spec R \ | \  x_1,\ldots,x_n\notin P\}$ and $\Spec R\mc{E}^{-1}$, so that since both extension and contraction send the zero ideal to the zero ideal, we get 
\begin{equation}\label{step a haon}
\ts{\bigcap_{P\in \Spec_d R,\ x_1,\ldots,x_n\notin P}P\neq 0\iff \bigcap_{Q\in \Spec_d R\mc{E}^{-1}}Q\neq 0.}
\end{equation}
We claim that 
\begin{equation}\label{step a do}
\ts{\bigcap_{P\in \Spec_d R,\ x_1,\ldots,x_n\notin P}P\neq 0\iff  \bigcap_{P\in \Spec_d R}P\neq 0.}
\end{equation}
One implication is trivial. For the other, suppose that $\bigcap_{P\in \Spec_d R,\ x_1,\ldots,x_n\notin P}P\neq 0$ and choose any $0\neq r$ which belongs to this intersection. Then $0\neq rx_1\cdots x_n\in \bigcap_{P\in \Spec_d R}P$, verifying \eqref{step a do}.
Now \eqref{step a haon} and \eqref{step a do} immediately give the result.
\end{proof}
\end{lemma}

\begin{lemma}\label{lemma localisation preserves weak gdme}
Let $R$ be a noetherian $\K$-algebra of finite Gelfand-Kirillov dimension in which every prime ideal is completely prime. If $R$ satisfies the quasi strong Dixmier-Moeglin equivalence and $\mc{E}$ is a right Ore set of regular elements of $R$ which is finitely generated as a multiplicative system, then $R\mc{E}^{-1}$ satisfies the quasi strong Dixmier-Moeglin equivalence. 
\begin{proof}
Every prime ideal of $R\mc{E}^{-1}$ takes the form $P\mc{E}^{-1}$ for some $P\in \Spec R$ with $P\cap \mc{E}=\emptyset$. Denoting by $\ol{\mc{E}}$ the image of $\mc{E}$ in $R/P$, we have 
\[
\begin{array}{rclr}
\locdeg P\mc{E}^{-1}&=& \locdeg \ang{0}_{R\mc{E}^{-1}/P\mc{E}^{-1}}& {} \\
                  {}&=& \locdeg \ang{0}_{(R/P)\ol{\mc{E}}^{-1}} & {}\\
                  {}&=& \locdeg \ang{0}_{R/P} & \tx{(Lemma \ref{localise and loc deg})}\\
                  {}&=& \locdeg P & {} \\
                  {}&=& \ratdeg P. & {} \\
\end{array}
\]
Since it is clear that $\ratdeg P = \ratdeg P\mc{E}^{-1}$, we are done.
\end{proof}
\end{lemma}

\begin{proposition}\label{transfer result for weak GDME}
Let $R$ be a noetherian $\K$-algebra of finite Gelfand-Kirillov dimension in which every prime ideal is completely prime. Suppose that for every $P\in \Spec R$, there exists a right Ore set $\mc{E}$ of regular elements of $R/P$ which is finitely generated as a multiplicative system, such that $(R/P)\mc{E}^{-1}$ satisfies the quasi strong Dixmier-Moeglin equivalence. Then $R$ itself satisfies the quasi strong Dixmier-Moeglin equivalence. 
\begin{proof}
Choose any $P\in \Spec R$. We have
\[
\begin{array}{rclr}
\locdeg P&=& \locdeg \ang{0}_{R/P} & {} \\
         &=& \locdeg \ang{0}_{(R/P)\mc{E}^{-1}} & \tx{(Lemma \ref{localise and loc deg})} \\
         &=& \ratdeg \ang{0}_{(R/P)\mc{E}^{-1}}. & {} \\
\end{array}
\]
Since it is clear that $\ratdeg \ang{0}_{(R/P)\mc{E}^{-1}}=\ratdeg P$, we are done. 
\end{proof}
\end{proposition}
\begin{remark}
The result of Proposition \ref{transfer result for weak GDME} holds if, rather than assuming that $(R/P)\mc{E}^{-1}$ satisfies the quasi strong Dixmier-Moeglin equivalence, we simply assume that, in $(R/P)\mc{E}^{-1}$, we have $\locdeg \ang{0}=\ratdeg \ang{0}$. 
\end{remark}

\section{The strong Dixmier-Moeglin equivalence in the commutative case}
\label{section:commutative}
If there is to be any hope that the strong Dixmier-Moeglin equivalence will hold for any quantum algebras, one should first check that it holds for commutative affine domains. Before checking this, let us introduce the useful notion of Tauvel's height formula:
\begin{definition}\label{THF}
\emph{Tauvel's height formula} is said to hold in a $\K$-algebra $R$ if for every prime ideal $P$ of $R$, the following equality holds:
\[
\GKdim R/P=\GKdim R-\height P.
\] 
\end{definition}
It is well known that Tauvel's height formula holds in commutative affine domains; as we shall remark later, it has also been shown to hold in several interesting quantum algebras, including all of those which interest us in this paper.
\begin{proposition}\label{GDME for affine domains}
Every commutative affine domain over $\K$ satisfies the strong Dixmier-Moeglin equivalence. 
\begin{proof}
Let $R$ be a commutative affine domain over $\K$ and let $P\in \Spec R$. We claim that 
\begin{equation}\label{kd}
\primdeg P=\kdim R/P=\ratdeg P.
\end{equation} 
Indeed, $R/P$ is itself a commutative affine domain, so that every primitive (i.e. maximal) ideal of $R/P$ has height $\kdim R/P$. It follows that $\primdeg P=\kdim R/P$. Moreover, by standard results of commutative algebra, we have
\[
\begin{array}{rcl}
                       \ratdeg P&=&\ts{\trdeg_\K}\mc{Z}(\Frac R/P)\\
                              {}&=&\ts{\trdeg_\K}\Frac(R/P)\\
                              {}&=&\kdim R/P,\\
\end{array}
\]
so that \eqref{kd} is proved.

Every maximal ideal of $R/P$ has height $\kdim R/P=\ratdeg P$, so that $\Spec_{1+\ratdeg P}R/P$ is empty and hence 
\[
\ts{\bigcap_{Q\in \Spec_{1+\ratdeg P}R/P}Q=R/P\neq 0.}
\] 
Now $\primdeg P=\ratdeg P\geq \locdeg P$. It is well known that commutative affine domains over $\K$ are catenary and satisfy the noncommutative Nullstellensatz over $\K$, so that by Proposition \ref{general ineq}, we have $\locdeg P\geq \primdeg P$. This completes the proof.
\end{proof}
\end{proposition}
\begin{remark}
Affine prime noetherian polynomial identity algebras over $\K$ can be shown to satisfy the strong Dixmier-Moeglin equivalence by a proof essentially the same as the proof above.
\end{remark}
\begin{remark}\label{kdim-ht affine domains}
Let $P$ be a prime ideal of a commutative affine domain $R$ over $\K$. Since Gelfand-Kirillov dimension and Krull dimension agree in commutative affine domains, Tauvel's height formula gives $\kdim R/P=\kdim R-\height P$. Now we conclude from Proposition \ref{GDME for affine domains} and equation \eqref{kd} that
\[
\locdeg P=\primdeg P=\ratdeg P = \kdim R-\height P.
\]
\end{remark}

\section{The strong Dixmier-Moeglin equivalence for uniparameter quantum tori}
\label{section:q-tori}
Let $N$ be a positive integer and let $\Lambda=(\lambda_{i,j})\in\mc{M}_N(\K^\times)$ be a multiplicatively skew-symmetric matrix. The \emph{quantum torus} associated to $\Lambda$ is denoted by $\mc{O}_{\Lambda}((\K^\times)^N)$ or $\K_\Lambda[T_1^{\pm 1},\ldots,T_N^{\pm 1}]$ and is presented as the $\K$-algebra generated by $T_1^{\pm 1},\ldots,T_N^{\pm 1}$ with relations 
\[
T_iT_i^{-1}=T_i^{-1}T_i=1\tx{ for all }i,\ T_jT_i=\lambda_{j,i}T_iT_j\tx{ for all }i,j.
\]
The algebra $\mc{O}_{\Lambda}((\K^\times)^N)$ can be written as the iterated skew-Laurent extension 
\[
\K[T_1^{\pm 1}][T_2^{\pm 1};\sigma_2]\cdots [T_N^{\pm 1};\sigma_N],
\]
where for each $j\in \llb 2,N \rrb$, $\sigma_j$ is the automorphism of $\K[T_1^{\pm 1}][T_2^{\pm 1};\sigma_2]\cdots [T_{j-1}^{\pm 1};\sigma_{j-1}]$ defined by $\sigma_j(T_i)=\lambda_{j,i}T_i$ for all $i\in \llb 1,j-1 \rrb$.
As such, $\mc{O}_{\Lambda}((\K^\times)^N)$ is a noetherian domain and there is a monomial $\K$-basis for $\mc{O}_{\Lambda}((\K^\times)^N)$ given by $\{T_1^{i_1}\cdots T_N^{i_N}\ | \ (i_1,\ldots,i_N)\in \Z^N\}$. By \cite[Corollary II.7.18]{BG}, $\mc{O}_{\Lambda}((\K^\times)^N)$ satisfies the noncommutative Nullstellensatz over $\K$ and by \cite[Theorem II.9.14]{BG}, $\mc{O}_{\Lambda}((\K^\times)^N)$ is catenary and satisfies Tauvel's height formula.

\vspace{2mm}

We recall from \cite[Section 1]{GL} some useful facts about quantum tori. 
For $\underline{i}=(i_1,\ldots,i_N)\in \Z^N$, we set $\underline{T}^{\underline{i}}:=T_1^{i_1}\cdots T_N^{i_N}$. For any $\ul{s},\ul{t}\in \Z^N$, we have $\ul{T}^{\ul{s}}\ul{T}^{\ul{t}}=\sigma(\ul{s},\ul{t})\ul{T}^{\ul{t}}\ul{T}^{\ul{s}}$, where $\sigma\co \Z^N\times \Z^N \to \K^\times$ is the alternating bicharacter which sends any $((s_1,\ldots,s_N),(t_1,\ldots,t_N))$ to $\prod_{i,j=1}^N\lambda_{i,j}^{s_it_j}$.

When $S$ is the subgroup $\{\ul{s}\in \Z^N \ | \ \sigma(\ul{s},-)\equiv 1\}$ of $\Z^N$, the centre of $\mc{O}_{\Lambda}((\K^\times)^N)$ is spanned over $\K$ by those $\ul{T}^{\ul{s}}$ with $\ul{s}\in S$. Where $\ul{b_1},\ldots,\ul{b_r}$ is a basis for $S$, the centre of $\mc{O}_{\Lambda}((\K^\times)^N)$ is a commutative Laurent polynomial ring in $(\ul{T}^{\ul{b_1}})^{\pm 1},\ldots,(\ul{T}^{\ul{b_r}})^{\pm 1}$.  Moreover, $\mc{O}_{\Lambda}((\K^\times)^N)$ is a free module over its centre with basis $\ul{T}^{\ul{t}}$, where $\ul{t}$ runs over any transversal for $S$ in $\Z^N$.

There is a bi-increasing homeomorphism, known as \emph{extension}, from $\Spec \mc{Z}(\mc{O}_{\Lambda}((\K^\times)^N))$ to $\Spec \mc{O}_{\Lambda}((\K^\times)^N)$ given by $I\mapsto \ang{I}$ (where $\ang{I}$ denotes the ideal of $\mc{O}_{\Lambda}((\K^\times)^N)$ generated by $I$). The inverse of this map is given by $J\mapsto J\cap \mc{Z}(\mc{O}_{\Lambda}((\K^\times)^N))$ and is known as \emph{contraction} from $\Spec \mc{O}_{\Lambda}((\K^\times)^N)$ to $\Spec \mc{Z}(\mc{O}_{\Lambda}((\K^\times)^N))$. In fact, contraction and extension define mutually inverse increasing bijections between the set of all ideals of $\mc{O}_{\Lambda}((\K^\times)^N)$ and the set of all ideals of its centre. 

\vspace{2mm}

Computing the rationality degree of a prime ideal $P$ of $\mc{O}_\Lambda((\K^\times)^N)$ requires study of the centre of $\Frac (\mc{O}_\Lambda((\K^\times)^N)/P)$. 
 The following general lemma is folklore, but we haven't been able to locate it in the literature\footnote{We thank Ken Goodearl for bringing this result to our attention.}. 
\begin{lemma}\label{Frac and Z commute}
Let $R$ be a noetherian domain and suppose that every nonzero ideal of $R$ intersects $\mc{Z}(R)$ nontrivially. Then 
\[
\mc{Z}(\Frac R)\cong \Frac\mc{Z}(R).
\] 
\begin{proof}
$\Frac\mc{Z}(R)$ embeds naturally into $\mc{Z}(\Frac R)$. Let $z\in \mc{Z}(\Frac R)$ and set $I=\{a\in R \ | \ za\in R\}$. Then $I$ is a nonzero ideal of $R$ and thus contains a nonzero element $c$ of $\mc{Z}(R)$. Now $z=(zc)c^{-1}\in \Frac\mc{Z}(R)$.  
\end{proof}
\end{lemma}
\begin{proposition}\label{Frac and Z commute for quots of quantum torus}
For a completely prime ideal $P$ of $\mc{O}_\Lambda((\K^\times)^N)$, we have
\[
\mc{Z}\left(\Frac \frac{ \mc{O}_\Lambda((\K^\times)^N)}{P}\right)\cong \Frac \mc{Z}\left(\frac{\mc{O}_\Lambda((\K^\times)^N)}{P}\right).
\] 
\begin{proof}
Set $R=\mc{O}_\Lambda((\K^\times)^N)$ and let $P$ be a completely prime ideal of $R$. 
By Lemma \ref{Frac and Z commute}, it will suffice to show that every nonzero ideal of $R/P$ intersects $\mc{Z}(R/P)$ nontrivially. This follows easily from the corresponding property of $R$ and the fact that every ideal of $R$ is generated by its intersection with $\mc{Z}(R)$.
\end{proof}
\end{proposition}

\begin{proposition}\label{quot centre}
For any ideal $I$ of $\mc{O}_\Lambda((\K^\times)^N)$, we have 
\[
\mc{Z}\left(\frac{\mc{O}_\Lambda((\K^\times)^N)}{I}\right)\cong \frac{\mc{Z}(\mc{O}_\Lambda((\K^\times)^N))}{I\cap \mc{Z}(\mc{O}_\Lambda((\K^\times)^N))}.
\]
\begin{proof}
Set $R=\mc{O}_\Lambda((\K^\times)^N)$. We may clearly assume that $I$ is a proper ideal and that $R$ is noncommutative. It follows that the subgroup $S$ of $\Z^N$ described at the beginning of this section has index at least three. Indeed $S$ clearly has index at least two (since $R$ is noncommutative), but if $S$ has index exactly two, then we may choose a transversal $\ul{0},\ul{i}$ for $S$ in $\Z^N$. Now $1=\ul{T}^{\ul{0}}$ and $\ul{T}^{\ul{i}}$ form a basis for $R$ as a module over its centre, from which it easily follows that $R$ is commutative, contradicting our assumption.

We claim that $\mc{Z}(R/I)=(\mc{Z}(R)+I)/I$. Indeed the inclusion $\mc{Z}(R/I)\supseteq (\mc{Z}(R)+I)/I$ is obvious. Suppose that $x\in R$ is central modulo $I$. 
We may choose elements $\ul{0},\ul{i_1},\ldots,\ul{i_n}$ ($n\geq 2$) of a transversal for $S$ in $\Z^N$ and central elements $z_0,z_1,\ldots,z_n$ of $R$ such that 
\[
x=z_0+\sum_{a=1}^nz_a\ul{T}^{\ul{i_a}}.
\]
We must have $\ul{T}^{\ul{i_1}}x(\ul{T}^{\ul{i_1}})^{-1}=x$ modulo $I$, so that 
\[
\sum_{a=1}^n(1-\sigma(\ul{i_1},\ul{i_a}))z_a\ul{T}^{\ul{i_a}}\in I
\]
and hence, by \cite[Proposition 1.4]{GL}, each $(1-\sigma(\ul{i_1},\ul{i_a}))z_a$ must belong to $I$. But since for $a\neq 1$ we have $\sigma(\ul{i_1},\ul{i_a})\neq 1$, we must have $z_2,\ldots,z_n\in I$ and hence $x=z_0+z_1\ul{T}^{\ul{i_1}}$ modulo $I$. 

We have $\ul{T}^{\ul{i_2}}x(\ul{T}^{\ul{i_2}})^{-1}=x$ modulo $I$ and using a similar argument to that which we employed above, we get $z_1\in I$ and hence $x=z_0$ modulo $I$, completing the proof. 
\end{proof} 
\end{proposition}
The quantum torus $\mc{O}_{\Lambda}((\K^\times)^N)$ is called a \emph{uniparameter} quantum torus if there exists a non root of unity $q\in \K^\times$ and an additively skew-symmetric matrix $A=(a_{i,j})\in \mc{M}_N(\Z)$ such that $\Lambda=(q^{a_{i,j}})$; in this case, we write $\mc{O}_{q,A}((\K^\times)^N)$ for $\mc{O}_\Lambda((\K^\times)^N)$. By \cite[Corollary II.6.10]{BG}, all prime ideals of $\mc{O}_{q,A}((\K^\times)^N)$ are completely prime so that Proposition \ref{Frac and Z commute for quots of quantum torus} applies. We are now ready to prove the main result of this section.
\begin{theorem}\label{GDME for upqt}
The uniparameter quantum tori $\mc{O}_{q,A}((\K^\times)^N)$ satisfy the strong Dixmier-Moeglin equivalence. 
\begin{proof}
Set $R=\mc{O}_{q,A}((\K^\times)^N)$ and choose any $P\in \Spec R$. As we have just noted, $P$ is guaranteed to be completely prime. Recall that $\mc{Z}(R)$ is a commutative Laurent polynomial ring; in particular, $\mc{Z}(R)$ is a commutative affine domain, so that it satisfies the strong Dixmier-Moeglin equivalence by Proposition \ref{GDME for affine domains}. By Propositions \ref{Frac and Z commute for quots of quantum torus} and \ref{quot centre}, we have
\[
\mc{Z}(\Frac R/P)\cong\Frac\mc{Z}(R/P)\cong \Frac\frac{\mc{Z}(R)}{\mc{Z}(R)\cap P}.
\]
It follows that $\ratdeg P=\ratdeg (\mc{Z}(R)\cap P)$. Since $\mc{Z}(R)/(\mc{Z}(R)\cap P)$ is a commutative affine domain, Remark \ref{kdim-ht affine domains} gives $\ratdeg P= \kdim \mc{Z}(R)-\height (\mc{Z}(R)\cap P)$.
Since extension and contraction are mutually inverse increasing homeomorphisms between $\Spec \mc{Z}(R)$ and $\Spec R$, we have $\height (\mc{Z}(R)\cap P)=\height P$, so that 
\[
\ratdeg P=\kdim \mc{Z}(R)-\height P.
\]
Every maximal ideal of $\mc{Z}(R)$ has height $\kdim \mc{Z}(R)$ and hence so does every maximal ideal of $R$. By \cite[Corollary 1.5]{GL}, the primitive ideals of $R$ are exactly its maximal ideals, so that every primitive ideal of $R$ has height $\kdim \mc{Z}(R)$. Now the catenarity of $R$ gives $\primdeg P=\kdim \mc{Z}(R)-\height P$ and, in particular, $\primdeg P=\ratdeg P$.

Since every maximal ideal of $R$ has height $\kdim \mc{Z}(R)$, every maximal ideal of $R/P$ has height $\kdim \mc{Z}(R)-\height P=\ratdeg P$. But then $\Spec_{1+\ratdeg P}R/P$ is empty so that\\ $\bigcap_{Q\in \Spec_{1+\ratdeg P}R/P}=R/P\neq 0$. This shows that $\locdeg P\leq \ratdeg P$. 

So far, we have shown that $\locdeg P\leq \primdeg P=\ratdeg P=\kdim \mc{Z}(R)-\height P$. Proposition \ref{general ineq} gives 
\begin{equation*}\label{all equal to kdim-ht}
\locdeg P=\primdeg P=\ratdeg P=\kdim \mc{Z}(R)-\height P. 
\end{equation*} 
\end{proof}
\end{theorem}

\section{Primer on $H$-stratification}
\label{section:H-stratification}
Our next aim is to show that uniparameter quantum affine spaces (which we shall later define) satisfy the strong Dixmier-Moeglin equivalence. For this, we will make use of the $H$-stratification theory of Goodearl and Letzter (for details on this theory, see \cite[II.2]{BG}). Indeed, an examination of the $H$-stratification (a notion which we define in this section) of a uniparameter quantum affine space reveals that every (prime homomorphic image of a) uniparameter quantum affine space localises to a (prime homomorphic image of a) uniparameter quantum torus. This allows us to transfer the quasi strong Dixmier-Moeglin equivalence from uniparameter quantum tori to uniparameter quantum affine spaces in Section \ref{section:q-affine}. Further examination of the $H$-stratification of a uniparameter quantum affine space allows us to calculate the primitivity degrees of the prime ideals and hence, in the next section, complete the proof that uniparameter quantum affine spaces satisfy the strong Dixmier-Moeglin equivalence.

The material in this section shall be useful beyond quantum affine spaces, so we work in a more general setting. Let us suppose that $R$ is a noetherian $\K$-algebra and that $H=(\K^\times)^r$ is an algebraic $\K$-torus acting rationally on $R$ by $\K$-algebra automorphisms. We refer to $H$-invariant prime ideals as \emph{$H$-prime} ideals. 
We denote by $H$-$\Spec R$ the \emph{$H$-spectrum} of $R$, namely the subspace of $\Spec R$ consisting of all $H$-prime ideals. Let us assume further that every $H$-prime ideal $J$ of $R$ is \emph{strongly $H$-rational} in the sense that the fixed field $\mc{Z}(\Frac(R/J))^H$ is $\K$ (in all of the algebras which will concern us in this paper, \cite[Theorem II.6.4]{BG} guarantees that every $H$-prime ideal is strongly $H$-rational).

For an ideal $I$ of $R$, $(I\co H):=\bigcap_{h\in H}h\cdot I$ is the largest $H$-invariant ideal of $R$ contained in $I$. It is well known that if $P$ is a prime ideal of $R$, then $(P\co H)$ is an $H$-prime ideal of $R$. 
For an $H$-prime ideal $J$ of $R$, the \emph{$H$-stratum} of $\Spec R$ associated to $J$ is denoted by $\Spec_JR$ and is defined by 
$\Spec_{J}R=\{P\in \Spec R \ | \ (P\co H)=J\}$.    
The $H$-strata form a partition of $\Spec R$, usually referred to as the $H$-\emph{stratification}. This stratification plays a crucial role in understanding the prime ideals of $R$ and, as we shall see later in this section, the primitive ideals of $R$. 
By \cite[Theorem II.2.13]{BG}, for each $H$-prime ideal $J$ of $R$, there is a bi-increasing homeomorphism from $\Spec_JR$ to the prime spectrum of an appropriate commutative Laurent polynomial algebra over $\K$; the \emph{Krull dimension} of the $H$-stratum $\Spec_J R$ is defined to be the Krull dimension of this commutative Laurent polynomial algebra.

Let us make a useful observation on the Krull dimension of $H$-strata under localisation. 
Let $\mc{E}$ be a right Ore set in $R$ consisting of regular $H$-eigenvectors with rational $H$-eigenvalues. There is a natural induced rational action of $H$ on $R\mc{E}^{-1}$ by $\K$-algebra automorphisms. 
Extension and contraction  restrict to mutually inverse increasing homeomorphisms between the set of $H$-prime ideals of $R$ which do not intersect $\mc{E}$ and the set of $H$-prime ideals of $R\mc{E}^{-1}$. Moreover, for any $H$-prime ideal $J$ of $R$ which does not intersect $\mc{E}$, extension and contraction restrict to mutually inverse increasing homeomorphisms between $\Spec_JR$ and $\Spec_{J\mc{E}^{-1}}R\mc{E}^{-1}$. 
We deduce:
\begin{lemma}\label{krull eq}
Let an algebraic $\K$-torus $H$ act rationally on a noetherian $\K$-algebra $R$ by $\K$-algebra automorphisms and suppose that all $H$-prime ideals of $R$ are strongly $H$-rational. Let $\mc{E}$ be a right Ore set in $R$ consisting of regular $H$-eigenvectors with rational $H$-eigenvalues. Then for any $H$-prime ideal $J$ of $R$ which does not intersect $\mc{E}$, we have
\[
\ts{\kdim \Spec_JR = \kdim \Spec_{J\mc{E}^{-1}}R\mc{E}^{-1}}.
\]
\end{lemma}

Under the further assumptions that $R$ has finitely many $H$-prime ideals and that $R$ satisfies the noncommutative Nullstellensatz over $\K$, \cite[Theorem II.8.4]{BG} says that $R$ satisfies the Dixmier-Moeglin equivalence and that the primitive ideals of $R$ are exactly those prime ideals which are maximal in their $H$-strata. Assuming further that $R$ is catenary and that the $H$-strata of $R$ satisfy a technical condition (given in inequality \eqref{the crucial inequality}), we now show that if $P$ is a prime ideal of $R$ belonging to $\Spec_JR$ for an $H$-prime ideal $J$ of $R$ and if $M\supseteq P$ is a primitive (i.e. maximal) element of $\Spec_JR$, then $\height M/P=\primdeg P$ (and we compute these quantities in terms of the Krull dimension of $\Spec_JR$). Crucially, this allows us to look only at a single $H$-stratum of $R$ in order to compute $\primdeg P$. 

\begin{proposition}\label{if we have the crucial inequality, then we know the prim deg}
Let $R$ be a catenary noetherian $\K$-algebra satisfying the noncommutative Nullstellensatz over $\K$ and let $H$ be an an algebraic $\K$-torus acting rationally on $R$ by $\K$-algebra automorphisms. Suppose that $H$-$\Spec R$ is finite, that all $H$-prime ideals of $R$ are strongly $H$-rational, and that for any pair of $H$-prime ideals $J\subseteq J'$ of $R$, we have 
\begin{equation}\label{the crucial inequality}
\ts{\kdim \Spec_{J}R+\height J\leq \kdim \Spec_{J'}R+\height J'.}
\end{equation}
Then for any $H$-prime ideal $J$ of $R$, any $P\in \Spec_JR$, and any primitive element $M\supseteq P$ of $\Spec_JR$, we have
\begin{equation}\label{prim deg computed}
\ts{\primdeg P=\height M/P=\kdim \Spec_JR +\height J-\height P.}
\end{equation}
\begin{proof}

Let $M$ be a primitive element of $\Spec_JR$ which contains $P$. Then $M$ is maximal in $\Spec_JR$, so that $\height M/J=\kdim \Spec_JR$. It follows from the catenarity of $R$ that 
\begin{equation}\label{ht of m}
\ts{\height M/P=\kdim \Spec_JR+\height J-\height P.}
\end{equation} 

Every primitive ideal of $R/P$ corresponds to a primitive ideal of $R$ which contains $P$. Choose any such primitive ideal $N$ of $R$ and say $N$ belongs to $\Spec_{J'}R$ for an $H$-prime ideal $J'$ of $R$. It is clear that $J\subseteq J'$.

Since $N$ is maximal in $\Spec_{J'}R$, we have $\height N/J'=\kdim \Spec_{J'}R$. It follows from the catenarity of $R$ that 
\begin{equation}\label{ht of n}
\ts{\height N/P=\kdim \Spec_{J'}R+\height J'-\height P.} 
\end{equation}
Equations \eqref{ht of m} and \eqref{ht of n}, along with the assumption \eqref{the crucial inequality}, show that the height of an arbitrary primitive ideal of $R/P$ is at least $\height M/P$. Since $M/P$ is itself primitive, we get $\height M/P=\primdeg P$; combining this with equation \eqref{ht of m} gives the result.
\end{proof}
\end{proposition}
\begin{remark}
Except for the inequality \eqref{the crucial inequality}, the conditions of Proposition \ref{if we have the crucial inequality, then we know the prim deg} are known to hold for many interesting algebras. Much of the rest of this paper is concerned with verifying inequality \eqref{the crucial inequality} for uniparameter quantum affine spaces (Section \ref{section:q-affine}) and quantum Schubert cells (Section \ref{section:q-schubert}). Our proofs rely on knowledge of the dimensions of the $H$-strata \cite{BCL, BL} and on knowledge of the posets of $H$-prime ideals \cite{GeigerYakimov, GL, CauchonMeriaux}.
\end{remark}

\section{The strong Dixmier-Moeglin equivalence for uniparameter quantum affine spaces}\label{Uniparameter quantum affine spaces}
\label{section:q-affine}

In a further step towards proving the strong Dixmier-Moeglin equivalence for quantum Schubert cells, we prove it in this section for uniparameter quantum affine spaces.

\subsection{Quantum affine spaces}

Let $N$ be a positive integer and let $\Lambda=(\lambda_{i,j})\in \mathcal{M}_N(\K^\times)$ be a multiplicatively skew-symmetric matrix. The \emph{quantum affine space} associated to $\Lambda$ is denoted by $\mc{O}_\Lambda(\K^N)$ or $\K_\Lambda[T_1,\ldots,T_N]$ and is presented as the $\K$-algebra with generators $T_1,\ldots,T_N$ and relations 
\[
T_jT_i=\lambda_{j,i}T_iT_j\tx{ for all }i,j\in \llb 1,N \rrb.
\]
The algebra $\mc{O}_\Lambda(\K^N)$ can be written as the iterated skew-polynomial extension 
\[
\K[T_1][T_2;\sigma_2]\cdots [T_N;\sigma_N],
\]
where, for each $j\in \llb 2,N \rrb$, $\sigma_j$ is the automorphism of $\K[T_1][T_2;\sigma_2]\cdots [T_{j-1};\sigma_{j-1}]$ defined by $\sigma_j(T_i)=\lambda_{j,i}T_i$ for all $i\in \llb 1,j-1\rrb$. As such, $\mc{O}_\Lambda(\K^N)$ is a noetherian domain. By \cite[Corollary II.7.18]{BG}, $\mc{O}_\Lambda(\K^N)$ satisfies the noncommutative Nullstellensatz over $\K$ and by \cite[Theorem II.9.14]{BG}, $\mc{O}_{\Lambda}(\K^N)$ is catenary and satisfies Tauvel's height formula. 
\subsection{$H$-stratification of $\Spec \mc{O}_\Lambda(\K^N)$}\label{strata in q affine space}
The algebraic $\K$-torus $H=(\K^\times)^N$ acts rationally on $\mc{O}_\Lambda(\K^N)$ by $\K$-algebra automorphisms as follows:
\[
(a_1,\ldots,a_N)\cdot T_i=a_iT_i\tx{ for all }i\in \llb 1,N \rrb \tx{ and all }(a_1,\ldots,a_N)\in H.
\]
For a subset $\Delta$ of $\{1,\ldots,N\}$, let $K_\Delta$ be the ideal of $\mc{O}_\Lambda(\K^N)$ generated by those $T_i$ with $i\in \Delta$. The ideal $K_\Delta$ is clearly an $H$-invariant completely prime ideal of $\mc{O}_\Lambda(\K^N)$. Goodearl and Letzter have shown \cite[Proposition 2.11]{GL} that all $H$-prime ideals of $\mc{O}_\Lambda(\K^N)$ take this form, namely that $H$-$\Spec \mc{O}_\Lambda(\K^N)=\{K_\Delta \ | \ \Delta\subseteq \{1,\ldots,N\}\}$. For any $\Delta\subseteq \{1,\ldots,N\}$, the $H$-stratum of $\mc{O}_\Lambda(\K^N)$ associated to $K_\Delta$ (which shall be denoted by $\Spec_\Delta (\mc{O}_\Lambda(\K^N))$) is given by 
\[
\ts{\Spec_\Delta}(\mc{O}_\Lambda(\K^N))=\left\{P\in \Spec \mc{O}_\Lambda(\K^N)\ | \ P\cap \{T_i\ | \ i\in \llb 1,N \rrb\}=\{T_i\ | \ i\in \Delta\}\right\}.
\]

\subsection{Uniparameter quantum affine spaces}\label{upqas}

The quantum affine space $\mc{O}_{\Lambda}(\K^N)$ is called a \emph{uniparameter} quantum affine space if there exists a non root of unity $q\in \K^\times$ and an additively skew-symmetric matrix $A=(a_{i,j})\in \mc{M}_N(\Z)$ such that $\Lambda=(q^{a_{i,j}})$. In this case, we denote $\mc{O}_\Lambda(\K^N)=\K_\Lambda[T_1,\ldots,T_N]$ by $\mc{O}_{q,A}(\K^N)=\K_{q,A}[T_1,\ldots,T_N]$. By \cite[Corollary II.6.10]{BG}, every prime ideal of $\mc{O}_{q,A}(\K^N)$ is completely prime. 

We use a transfer result from Section \ref{section:generalities} to show that $\mc{O}_{q,A}(\K^N)$ satisfies the quasi strong Dixmier-Moeglin equivalence.
\begin{proposition}\label{weak gdme for qas}
The uniparameter quantum affine spaces $\mc{O}_{q,A}(\K^N)$ satisfy the quasi strong Dixmier-Moeglin equivalence.
\begin{proof}
Set $R=\mc{O}_{q,A}(\K^N)=\K_{q,A}[T_1,\ldots,T_N]$. Choose any $P\in \Spec R$ and say $P\in \Spec_\Delta R$ for a subset $\Delta$ of $\{1,\ldots,N\}$. Let $\mc{E}$ be the multiplicative system in $R$ generated by those $T_i$ for which $i\notin \Delta$. Then $\mc{E}$ satisfies the Ore condition on both sides in $R$ and, denoting by $\ol{\mc{E}}$ and $\hat{\mc{E}}$ its images in $R/P$ and $R/K_\Delta$ respectively, we have
\[
(R/P)\ol{\mc{E}}^{-1}\cong ((R/K_\Delta)\hat{\mc{E}}^{-1})/((P/K_\Delta)\hat{\mc{E}}^{-1}).
\]
The uniparameter quantum torus $(R/K_\Delta)\hat{\mc{E}}^{-1}$ satisfies the strong Dixmier-Moeglin equivalence by Theorem \ref{GDME for upqt} and hence so does its homomorphic image $(R/P)\ol{\mc{E}}^{-1}$. The result now follows from Proposition \ref{transfer result for weak GDME}.
\end{proof}
\end{proposition} 

\subsection{The strong Dixmier-Moeglin equivalence for uniparameter quantum affine spaces}
Since we have proven that $\mc{O}_{q,A}(\K^N)$ satisfies the quasi strong Dixmier-Moeglin equivalence, proving that $\primdeg P=\ratdeg P$ holds for all prime ideals $P$ of $\mc{O}_{q,A}(\K^N)$ will establish the strong Dixmier-Moeglin equivalence for $\mc{O}_{q,A}(\K^N)$. 

In order to invoke Proposition \ref{if we have the crucial inequality, then we know the prim deg}, which gives us an expression for the primitivity degree of any prime ideal $P$ of $\mc{O}_{q,A}(\K^N)$ in terms of the dimension of the $H$-stratum to which $P$ belongs, we must prove an inequality relating the dimensions of $H$-strata of $\mc{O}_{q,A}(\K^N)$. First we introduce some new notation:
\begin{notation}
Let $\Delta$ be a  subset of $\{1,\ldots,N\}$ and set $\{\ell_1<\ldots<\ell_d\}=\{1,\ldots,N\} \bs \Delta$. We define the \emph{skew-adjacency matrix}, $A(\Delta)$, of $\Delta$ to be the $d\times d$ additively skew-symmetric submatrix of $A=(a_{i,j})\in \mc{M}_N(\Z)$ whose $(s,t)$ entry $(s<t)$ is $a_{\ell_s,\ell_t}$.
\end{notation}
For any subset $\Delta$ of $\{1,\ldots,N\}$, it follows from \cite[Theorem 3.1]{BL} that the dimension of the $H$-stratum $\Spec_\Delta (\mc{O}_{q,A}(\K^N))$ corresponding to the $H$-prime ideal $K_\Delta=\ang{T_i\ | \ i\in \Delta}$ is exactly $\dim_\Q (\ker A(\Delta))$. In fact, \cite[Theorem 3.1]{BL} applies to a more general class of algebras called \emph{uniparameter Cauchon-Goodearl-Letzter extensions} (see Section \ref{CGL DDA}).

\begin{proposition}\label{crucial inequality for qas}
For any pair of $H$-prime ideals $K_\Delta\subseteq K_{\Delta'}$ of $\mc{O}_{q,A}(\K^N)$, we have 
\[
\ts{\kdim \Spec_\Delta (\mc{O}_{q,A}(\K^N))+\height K_\Delta\leq \kdim \Spec_{\Delta'} (\mc{O}_{q,A}(\K^N))+\height K_{\Delta'}.}
\]
\begin{proof}
Since $K_\Delta\subseteq K_{\Delta'}$, we clearly have $\Delta\subseteq \Delta'$. The matrix $A(\Delta')$ is an $(N-|\Delta'|)$-square submatrix of the $(N-|\Delta|)$-square matrix $A(\Delta)$, so that $\rk A(\Delta')\leq \rk A(\Delta)$ and 
\[
(N-|\Delta'|)-\dim_\Q(\ker A(\Delta'))\leq (N-|\Delta|)-\dim_\Q(\ker A(\Delta)).
\]
Hence, we have 
\begin{equation}\label{card}
\dim_\Q(\ker A(\Delta))+|\Delta|\leq \dim_\Q(\ker A(\Delta'))+|\Delta'|.
\end{equation}
Tauvel's height formula holds in $\mc{O}_{q,A}(\K^N)$, so that 
\[
\height K_\Delta=\GKdim \mc{O}_{q,A}(\K^N)-\GKdim (\mc{O}_{q,A}(\K^N)/K_\Delta)=N-(N-|\Delta|)=|\Delta|
\] and similarly $\height K_{\Delta'}=|\Delta'|$. Now \eqref{card} and \cite[Theorem 3.1]{BL} give 
\[
\ts{\kdim \Spec_\Delta (\mc{O}_{q,A}(\K^N))+\height K_\Delta\leq \kdim \Spec_{\Delta'} (\mc{O}_{q,A}(\K^N))+\height K_{\Delta'}.}
\]
\end{proof} 
\end{proposition}

With Proposition \ref{crucial inequality for qas} in hand, we can apply Proposition \ref{if we have the crucial inequality, then we know the prim deg} to $\mc{O}_{q,A}(\K^N)$ in our proof of the main result of this section:

\begin{theorem}\label{gdme for qas}
The uniparameter quantum affine spaces $\mc{O}_{q,A}(\K^N)$ satisfy the strong Dixmier-Moeglin equivalence.
\begin{proof}
Set $R=\mc{O}_{q,A}(\K^N)=\K_{q,A}[T_1,\ldots,T_N]$. We showed in Proposition \ref{weak gdme for qas} that $R$ satisfies the quasi strong Dixmier-Moeglin equivalence, so what remains is to prove that $\primdeg P=\ratdeg P$ for all prime ideals $P$ of $R$.  

Let $P$ be any prime ideal of $R$ and say $P\in \Spec_\Delta R$ for a subset $\Delta$ of $\{1,\ldots,N\}$. Proposition \ref{if we have the crucial inequality, then we know the prim deg} gives \[
\ts{\primdeg P=\kdim \Spec_\Delta R+\height K_\Delta-\height P.} 
\]

Let $\mc{E}$ be the multiplicative system in $R$ generated by those $T_i$ for which $i\notin \Delta$. Then $\mc{E}$ satisfies the Ore condition on both sides in $R$ and, denoting by $\hat{\mc{E}}$ its image in $R/K_\Delta$, we have $R\mc{E}^{-1}/P\mc{E}^{-1}\cong ((R/K_\Delta)\hat{\mc{E}}^{-1})/((P/K_\Delta)\hat{\mc{E}}^{-1})$. Notice that $(R/K_\Delta)\hat{\mc{E}}^{-1}$ is a uniparameter quantum torus and that $P\mc{E}^{-1}\in \Spec_{K_\Delta \mc{E}^{-1}}R\mc{E}^{-1}$.

Since $R$ is catenary and noetherian, so is $R\mc{E}^{-1}$. Moreover, $R\mc{E}^{-1}$ can be obtained from $\K$ by a finite number of skew-polynomial and skew-Laurent extensions; in particular, $R\mc{E}^{-1}$ is a constructible $\K$-algebra in the sense of \cite[9.4.12]{McConnellRobson}, so that $R\mc{E}^{-1}$ satisfies the noncommutative Nullstellensatz over $\K$ by \cite[Theorem 9.4.21]{McConnellRobson}. From the discussion of the effect of localisation on $H$-stratification (Section \ref{section:H-stratification}), we deduce that $R\mc{E}^{-1}$ satisfies the conditions of Proposition \ref{if we have the crucial inequality, then we know the prim deg} and that \begin{align*}
\primdeg P\mc{E}^{-1}&=\ts{\kdim \Spec_{K_\Delta\mc{E}^{-1}} R\mc{E}^{-1}}+\height K_\Delta\mc{E}^{-1}-\height P\mc{E}^{-1} \\
                     &=\ts{\kdim \Spec_\Delta R}+\height K_\Delta-\height P, \\
\end{align*}
so that 
\[
\primdeg P\mc{E}^{-1}=\primdeg P.
\]

Since the uniparameter quantum torus $(R/K_\Delta)\hat{\mc{E}}^{-1}$ satisfies the strong Dixmier-Moeglin equivalence (Theorem \ref{GDME for upqt}), so does its homomorphic image $R\mc{E}^{-1}/P\mc{E}^{-1}$. So $\primdeg \ang{0}=\ratdeg \ang{0}$ holds in $R\mc{E}^{-1}/P\mc{E}^{-1}$, which can be rephrased by saying that in $R\mc{E}^{-1}$, we have $\primdeg P\mc{E}^{-1}=\ratdeg P\mc{E}^{-1}$. Since we have already shown that $\primdeg P=\primdeg P\mc{E}^{-1}$ and it is clear that $\ratdeg P\mc{E}^{-1}=\ratdeg P$, we have $\primdeg P=\ratdeg P$, as required.
\end{proof}
\end{theorem}

\section{CGL extensions and the deleting derivations algorithm}\label{CGL DDA}
In the terminology introduced in \cite[Definition 3.1]{LLR}, let $R=\K[X_1][X_2;\sigma_2,\delta_2]\cdots [X_N;\sigma_N;\delta_N]$ be a uniparameter \emph{Cauchon-Goodearl-Letzter (CGL)} extension. This class of algebras contains many quantum algebras such as quantum matrices and, more generally, quantum Schubert cells. In particular, there exists an algebraic $\K$-torus $H=(\K^\times)^d$ acting rationally on $R$ by $\K$-algebra automorphisms, there exists $q\in \K^\times$ not a root of unity, and there exists an additively skew-symmetric matrix $A=(a_{i,j})\in \mc{M}_N(\Z)$ such that 
 \begin{enumerate}[(i)]
\item For all $j\in \llb 2,N \rrb$, $\delta_j$ is locally nilpotent;
\item For all $j\in \llb 2,N \rrb$, there exists $q_j\in \K^\times$ not a root of unity such that $\sigma_j\circ \delta_j=q_j\delta_j\circ \sigma_j$;  
\item For all $j\in \llb 2,N \rrb$ and all $i\in \llb 1,j-1\rrb$, we have $\sigma_j(X_i)=q^{a_{j,i}}X_i$;
\item $X_1,\ldots,X_N$ are $H$-eigenvectors;
\item The set $\{\lambda\in \K^\times \ | \ \tx{there exists }h\in H\tx{ such that }h\cdot X_1=\lambda X_1\}$ is infinite;
\item For all $j\in \llb 2,N \rrb$, there exists $h_j\in H$ such that $h_j\cdot X_j=q_jX_j$ and, for $i\in \llb 1,j-1\rrb$, $h_j\cdot X_i=q^{a_{j,i}}X_i$.
\end{enumerate}
$R$ is a noetherian domain and it satisfies the noncommutative Nullstellensatz over $\K$ by \cite[Theorem II.7.17]{BG}. By \cite[Theorem II.6.9]{BG}, all prime ideals of $R$ are completely prime.

 Cauchon \cite{Cauchon} introduced an algorithm (now known as the \emph{deleting derivations algorithm}) which relates the prime spectrum and the $H$-stratification of $R$ to those of the quantum affine space $\ol{R}$ which results from ``deleting" the derivations $\delta_i$ (for a survey of this algorithm, see \cite[Section 2C]{BL}). Following the notation of \cite{Cauchon}, $\ol{R}$ is, more precisely, a uniparameter quantum affine space in indeterminates $T_1,\ldots,T_N$ with commutation relations given by $q$ and the matrix $A$, i.e. 
\[
\ol{R}=\K_{q,A}[T_1,\ldots,T_N]=\mc{O}_{q,A}(\K^N).
\]
There is a \emph{canonical injection} $\varphi$ of $\Spec R$ into $\Spec \ol{R}$ (see \cite[Section 4]{Cauchon}), which Cauchon used to construct a partition of $\Spec R$ which we now describe. 

Let $W$ be the power set of $\{1,\ldots,N\}$. For any $\Delta\in W$, set $\Spec_\Delta R=\varphi^{-1}(\Spec_{\Delta}\ol{R})$, where $\Spec_\Delta \ol{R}$ denotes the stratum in $\Spec \ol{R}$ associated to the $H$-prime ideal $K_\Delta=\ang{T_i\ | \ i\in \Delta}$ (see Subsection \ref{strata in q affine space}). Denote by $W'$ the set of those $\Delta\in W$ with $\Spec_\Delta R\neq \emptyset$. The elements of $W$ are called the \emph{diagrams} of the CGL extension $R$ and the elements of $W'$ are called the \emph{Cauchon diagrams} of $R$.
By \cite[Proposition 4.4.1]{Cauchon}, we have
\[
\Spec R=\bigsqcup_{\Delta\in W'}\ts{\Spec_\Delta R.}
\]
This is called the \emph{canonical partition} of $\Spec R$ and, by \cite[Th\'eor\`eme 5.5.2]{Cauchon}, it coincides with the partition of $\Spec R$ into $H$-strata. Let us make this more precise. 

For any Cauchon diagram $\Delta$ of $R$, the canonical injection $\varphi$ restricts to a bi-increasing homeomorphism from $\Spec_\Delta R$ to $\Spec_\Delta \ol{R}$ (\cite[Th\'eor\`emes 5.1.1 and 5.5.1]{Cauchon}). Moreover, by \cite[Lemme 5.5.8 and Th\'eor\`eme 5.5.2]{Cauchon},  we have the following description of the $H$-prime ideals of $R$:
\begin{enumerate}[(i)]
\item For any $\Delta\in W'$, there is a (unique) $H$-invariant (completely) prime ideal $J_\Delta$ of $R$ such that $\varphi(J_\Delta)=K_\Delta$;
\item $H$-$\Spec R=\{J_\Delta \ | \ \Delta\in W'\}$;
\item $\Spec_{J_\Delta}R=\Spec_\Delta R$ for all $\Delta\in W'$.
\end{enumerate}
The invertible map $\Delta\mapsto J_\Delta$ from $W'$ to $H$-$\Spec R$ is increasing but, in general, its inverse is not.

We are now in position to establish the quasi strong Dixmier-Moeglin equivalence for uniparamater CGL extensions.

\begin{theorem}\label{uniparameter CGLs satisfy the weak GDME}
Every uniparameter CGL extension satisfies the quasi strong Dixmier-Moeglin equivalence. 
\begin{proof}
Let $R$ be a  uniparameter CGL extension. Recall that both in $R$ and in the uniparameter quantum affine space $\ol{R}$, all prime ideals are completely prime.

Choose any $P\in \Spec R$ and say $P\in \Spec_\Delta R$ for a Cauchon diagram $\Delta$ of $R$. Let $\mc{E}$ be the image in $\ol{R}/\varphi(P)$ of the multiplicative system in $\ol{R}$ generated by those $T_i$ for which $i\in \{1,\ldots,N\} \bs \Delta$. By \cite[Th\'eorem\`e 5.4.1]{Cauchon}, $\mc{E}$ satisfies the Ore condition on both sides in $\ol{R}/\varphi(P)$ and there exists a finitely generated multiplicative system $\mc{F}$  
in $R/P$ satisfying the Ore condition on both sides such that 
\begin{equation}\label{an isom}
(R/P)\mc{F}^{-1}\cong (\ol{R}/\varphi(P))\mc{E}^{-1}.
\end{equation} 
 
Since $\ol{R}$ is a uniparameter quantum affine space, it satisfies the strong Dixmier-Moeglin equivalence (Theorem \ref{gdme for qas}) and hence so does every homomorphic image of $\ol{R}$. In particular, $\ol{R}/\varphi(P)$ satisfies the strong Dixmier-Moeglin equivalence. Hence, by Lemma \ref{lemma localisation preserves weak gdme}, $ (\ol{R}/\varphi(P))\mc{E}^{-1}$ satisfies the quasi strong Dixmier-Moeglin equivalence. 
The result now follows from \eqref{an isom} and Proposition \ref{transfer result for weak GDME}.
\end{proof}
\end{theorem}

Regarding the strong Dixmier-Moeglin equivalence, we can prove the following partial result.

\begin{theorem}\label{GDME for CGLs with crucial ineq}
 If $R$ is a catenary uniparameter CGL extension such that for any pair of $H$-prime ideals $J\subseteq J'$ of $R$, the following inequality holds:
\begin{equation}\label{assume we have the ineq}
\ts{\kdim \Spec_JR+\height J\leq \kdim \Spec_{J'}R+\height J'},
\end{equation}
then $R$ satisfies the strong Dixmier-Moeglin equivalence.
\begin{proof}
Since $R$ satisfies the quasi strong Dixmier-Moeglin equivalence (Theorem \ref{uniparameter CGLs satisfy the weak GDME}), we need only show that for every prime ideal $P$ of $R$, we have $\primdeg P=\ratdeg P$. 
By \cite[Theorem II.8.4]{BG}, $R$ and $\ol{R}$ satisfy the Dixmier-Moeglin equivalence and, in each of these two algebras, the primitive ideals are exactly the prime ideals which are maximal in their $H$-strata.

Suppose that $P$ is a prime ideal of $R$ with $P\in \Spec_\Delta R$ for a Cauchon diagram $\Delta$ of $R$. Choose any primitive (i.e. maximal) element $M\supseteq P$ of $\Spec_\Delta R$. Since $\varphi$ restricts to a bi-increasing homeomorphism from $\Spec_\Delta R$ to $\Spec_\Delta \ol{R}$, we get that $\varphi(M)$ is a maximal (i.e. primitive) element of $\Spec_\Delta \ol{R}$ and that $\varphi(M)$ contains $\varphi(P)$. Proposition \ref{crucial inequality for qas} and the assumption \eqref{assume we have the ineq} allow us to invoke Proposition \ref{if we have the crucial inequality, then we know the prim deg} to get  
\begin{equation}\label{using correct strata result}
 \height M/P=\primdeg P  \tx{ and }   \height \varphi(M)/\varphi(P)=\primdeg \varphi(P). 
\end{equation}

Moreover, since $\varphi$ restricts to a bi-increasing homeomorphism from $\Spec_\Delta R$ to $\Spec_\Delta \ol{R}$, it induces a length-preserving one-to-one correspondence between the chains of prime ideals from $P$ to $M$ and the chains of prime ideals from $\varphi(P)$ to $\varphi(M)$. It follows that 
\begin{equation}\label{varphi preserves chains}
\height M/P=\height \varphi(M)/\varphi(P).
\end{equation} 
We deduce from \eqref{using correct strata result} and \eqref{varphi preserves chains} that $\primdeg P=\primdeg \varphi(P)$. Now, recalling that the uniparameter quantum affine space $\ol{R}$ satisfies the strong Dixmier-Moeglin equivalence (by Theorem \ref{gdme for qas}) and that, by \cite[Th\'eorem\`e 5.4.1]{Cauchon}, $\Frac(R/P)\cong \Frac(\ol{R}/\varphi(P))$, we have 
\begin{align*}
\primdeg P&=\primdeg \varphi(P)\\
{ }&=\ratdeg \varphi(P)\\
{ }&=\ratdeg P,\\
\end{align*}
as required.
\end{proof}
\end{theorem}

\section{Quantum Schubert cells} 
\label{section:q-schubert}

We discuss quantum Schubert cells and their uniparameter CGL extension structure. Yakimov \cite[Theorem 5.7]{Yakimov} has shown that these algebras are catenary and satisfy Tauvel's height formula. We show that they satisfy inequality \eqref{assume we have the ineq} so that, by Theorem \ref{GDME for CGLs with crucial ineq}, they satisfy the strong Dixmier-Moeglin equivalence.
\subsection{The algebras $U_q[w]$ and their uniparameter CGL extension structure}

Let $\g$ be a simple complex Lie algebra of rank $n$ and let $\pi:=\{\alpha_1,\ldots,\alpha_n\}$ be the set of simple roots associated to a triangular decomposition $\g=\n^-\oplus \h\oplus \n^+$ of $\g$. The set $\pi$ is a basis of a real Euclidean vector space $E$, whose inner product we denote by $(-,-)$. Recall that the Weyl group of $\g$, which we denote by $\mc{W}$, is the subgroup of the orthogonal group of $E$ generated by the reflections $s_i:=s_{\alpha_i}$, for $i=1,\ldots,n$, with reflecting hyperplanes $H_i:=\{\beta\in E\ |\ (\beta,\alpha_i)=0\}$, $i=1,\ldots,n$.  

Where $q\in \K^\times$ is not a root of unity and $w$ is any element of $\mc{W}$, De Concini, Kac, and Procesi \cite{DKP} 
defined a quantum analogue, $U_q[w]$, of the universal enveloping algebra of the nilpotent Lie algebra $\n^+\cap \Ad_w(\n^-)$, where $\Ad$ denotes the adjoint action.
We refer the reader to \cite[Subsection 3C]{BL} for a description of the \emph{quantum Schubert cell} $U_q[w]$ as a certain subalgebra of $U_q^+(\g)$, where $U_q(\g)$ is the quantised enveloping algebra of $\g$ over $\K$ associated to the above data. 

 $\mc{W}$ is a Coxeter group with respect to the generators $s_1,\ldots,s_n$ and we define the length, $\ell(w)$, of $w$ to be the smallest $N$ such that there exist $i_j\in \{1,\ldots,n\}$ satisfying $w=s_{i_1}\cdots s_{i_N}$. Let us fix this \emph{reduced} expression for $w$.
It is well known that $\beta_1=\alpha_{i_1},\ \beta_2=s_{i_1}(\alpha_{i_2}),\ldots,\beta_N=s_{i_1}\cdots s_{i_{N-1}}(\alpha_{i_N})$ are distinct positive roots and that the set $\{\beta_1,\ldots,\beta_N\}$ does not depend on the chosen reduced expression for $w$.

Cauchon proved \cite[Proposition 6.1.2 and Lemme 6.2.1]{Cauchon} that $U_q[w]$ is a uniparameter CGL extension in $N$ indeterminates with the following associated additively skew-symmetric matrix:
\begin{equation}\label{Schubert matrix}
A:=\left(
\begin{array}{ccccc}
0&(\beta_1,\beta_2)&\cdots&\cdots&(\beta_1,\beta_N)\\
-(\beta_1,\beta_2)&0&(\beta_2,\beta_3)&{}&(\beta_2,\beta_N)\\
\vdots&\ddots&\ddots&\ddots&\vdots\\
\vdots&{}&\ddots&0&(\beta_{N-1},\beta_N)\\
-(\beta_1,\beta_N)&\cdots&\cdots&-(\beta_{N-1},\beta_N)&0\\
\end{array}
\right).
\end{equation}
Theorem \ref{uniparameter CGLs satisfy the weak GDME} immediately gives:
\begin{proposition}
The quantum Schubert cells $U_q[w]$ satisfy the quasi strong Dixmier-Moeglin equivalence. 
\end{proposition}

\subsection{The strong Dixmier-Moeglin equivalence for $U_q[w]$} 

Considering $U_q[w]$ as a uniparameter CGL extension in $N$ indeterminates with associated additively skew-symmetric matrix $A$ (see \eqref{Schubert matrix}), recall that $J_\Delta$ denotes the $H$-prime ideal of $U_q[w]$ associated to a Cauchon diagram $\Delta$ of $U_q[w]$. The remaining work lies in proving that for any pair of $H$-prime ideals $J_\Delta\subseteq J_{\Delta'}$ of $U_q[w]$, the following inequality holds:
\[
\ts{\kdim \Spec_{{\Delta}} U_q[w]+\height J_\Delta\leq \kdim \Spec_{{\Delta'}} U_q[w]+\height J_{\Delta'}.}
\] 
This will allow us to invoke Theorem \ref{GDME for CGLs with crucial ineq} to show that $U_q[w]$ satisfies the strong Dixmier-Moeglin equivalence. 

In contrast to that of most algebras supporting an $H$-action, the poset structure of the $H$-spectrum of $U_q[w]$ is known. Let us denote by $\leq$ the Bruhat order on $\mc{W}$ and let us set $\mc{W}^{\leq w}:=\{u\in \mc{W}\ | \ u\leq w\}$. The posets $H$-$\Spec U_q[w]$ and $\mc{W}^{\leq w}$ are isomorphic. In order to describe an isomorphism, we introduce some notation:
\begin{notation}
Recall that we have fixed a reduced expression $w=s_{i_1}\cdots s_{i_N}$ for $w$. Let $\Delta \subseteq \{1,\ldots,N\}$ be any (not necessarily Cauchon) diagram.
\begin{enumerate}[(i)]
\item For all $k=1,\ldots,N$, we set 
\[
\piecewise{s_{i_k}^\Delta :}{s_{i_k}}{\tx{if }k\in \Delta}{\id}{\tx{otherwise}.}
\]
\item We set $\{l_1<\cdots <l_d\}:=\{1,\ldots,N\}\bs \Delta$ and $j_r=i_{l_r}$ for all $r=1,\ldots,d$.
\item We set $w^\Delta:=s_{i_1}^\Delta\cdots s_{i_N}^\Delta\in \mc{W}$.
\item We set $A(w^\Delta)$ to be the $d\times d$ additively skew-symmetric submatrix of $A$ whose $(s,t)$-entry ($s<t$) is $(\beta_{j_s},\beta_{j_t})$.
\end{enumerate}
\end{notation}

Cauchon and M\'eriaux \cite[Corollary 5.3.1]{CauchonMeriaux} showed that the map
\begin{equation}\label{pos iso}
\tx{$H$-$\Spec U_q[w]$}\to \mc{W}^{\leq w};\ J_\Delta\mapsto w^\Delta,
\end{equation}
where $\Delta$ runs over the set of Cauchon diagrams of $U_q[w]$, is a bijection; they asked whether or not this bijection is an isomorphism of posets and this question was answered affirmatively by Geiger and Yakimov \cite[Theorem 4.4]{GeigerYakimov}. 
  
\begin{lemma}\label{ht of H primes}
For any Cauchon diagram $\Delta$ of $U_q[w]$, we have $\height J_\Delta=|\Delta|$.
\begin{proof}
Set $R=U_q[w]$ and recall that $\ol{R}$ denotes the uniparameter quantum affine space (in indeterminates $T_1,\ldots,T_N$ say) which results from ``deleting" the derivations in the expression of $R$ as a uniparameter CGL extension in $N$ indeterminates. Recall that $K_\Delta=\ang{T_i\ |\ i\in \Delta}$ is the image of $J_\Delta$ under the canonical injection $\varphi\co \Spec R\to \Spec \ol{R}$. 

Let $\mc{E}$ be the image in $\ol{R}/K_\Delta$ of the multiplicative system in $\ol{R}$ generated by those $T_i$ for which $i\notin \Delta$. Then $\mc{E}$ satisfies the Ore condition on both sides in $\ol{R}/K_\Delta$ and it follows from \cite[Th\'eorem\`e 5.4.1]{Cauchon} both that $R/J_\Delta$ embeds in the uniparameter quantum torus $(\ol{R}/K_\Delta)\mc{E}^{-1}$  and that $\Frac (R/J_\Delta)\cong \Frac ((\ol{R}/K_\Delta)\mc{E}^{-1})$. By \cite[Corollary 2.2]{LorenzTDeg}, the uniparameter quantum torus $(\ol{R}/K_\Delta)\mc{E}^{-1}$ is Tdeg-stable (in the sense of \cite[Section 1]{Zhang}). Therefore, we can apply \cite[Proposition 3.5(4)]{Zhang} to get $\GKdim R/J_\Delta=\GKdim (\ol{R}/K_\Delta)\mc{E}^{-1}=N-|\Delta|$.

Since $R$ satisfies Tauvel's height formula, we conclude that
$$N-|\Delta |=\GKdim R/J_\Delta=\GKdim R-\height J_\Delta=N-\height J_\Delta,$$
and so $\height J_\Delta=|\Delta|$, as desired. 
\end{proof}
\end{lemma}

We are now in position to establish the crucial inequality required to prove that quantum Schubert cells satisfy the strong Dixmier-Moeglin equivalence.

\begin{proposition}\label{crucial inequality for qsc}
For any pair of $H$-prime ideals $J_\Delta\subseteq J_{\Delta'}$ of $U_q[w]$, we have 
\[
\ts{\kdim \Spec_{{\Delta}} U_q[w]+\height J_\Delta\leq \kdim \Spec_{{\Delta'}} U_q[w]+\height J_{\Delta'}.}
\]
\begin{proof}
As we have noted, $U_q[w]$ is a uniparameter CGL extension in $N$ indeterminates with associated additively skew-symmetric matrix $A$. 
By \cite[Theorems 2.3 and 3.1]{BCL}, we have
\[
\ts{\kdim \Spec_{\Delta}} U_q[w]=\dim_\Q \ker (w^\Delta+w)\tx{ and } \ts{\kdim \Spec_{{\Delta'}}}U_q[w]=\dim_\Q \ker A(w^{\Delta'}).
\]
From the poset isomorphism ($H$-$\Spec U_q[w]\to \mc{W}^{\leq w};\ J_\Delta\mapsto w^\Delta$), we deduce that $w^\Delta\leq w^{\Delta'}$. Since the diagrams $\Delta$ and $\Delta'$ are Cauchon, the subexpressions $w^\Delta$ and $w^{\Delta'}$ of $w=s_{i_1}\cdots s_{i_N}$ are reduced by \cite[Corollary 5.3.1(2)]{CauchonMeriaux}. Since $w^\Delta\leq w^{\Delta'}$, \cite[Corollary 5.8]{HumphreysCoxeterGroups} allows us to choose a diagram (not necessarily Cauchon) $\widetilde{\Delta}\subseteq \Delta'$ such that $w^{\widetilde{\Delta}}=w^\Delta$ and the subexpression $w^{\widetilde{\Delta}}$ of $w=s_{i_1}\cdots s_{i_N}$ is reduced. Now $\ts{\kdim \Spec_{{\Delta}}} U_q[w]=\dim_\Q \ker (w^{\widetilde{\Delta}}+w)$, 
so that \cite[Theorem 3.1]{BCL} gives $\ts{\kdim \Spec_{{\Delta}}} U_q[w]=\dim_\Q \ker A(w^{\widetilde{\Delta}})$.

$A(w^{\Delta'})$ is an $(N-|\Delta'|)$-square submatrix of the $(N-|\widetilde{\Delta}|)$-square matrix $A(w^{\widetilde{\Delta}})$, so that 
$\rk A(w^{\Delta'}) \leq \rk A(w^{\widetilde{\Delta}})$ and hence $\dim_\Q \ker A(w^{\widetilde{\Delta}})+|\widetilde{\Delta}|\leq \dim_\Q \ker A(w^{\Delta'})+|\Delta'|$ and
\begin{equation}\label{almost}
\ts{\kdim \Spec_{{\Delta}}} U_q[w]+|\widetilde{\Delta}|\leq \ts{\kdim \Spec_{{\Delta'}}} U_q[w]+|\Delta'|.
\end{equation}
 By Lemma \ref{ht of H primes}, we have $\height J_\Delta=|\Delta|$ and $\height J_{\Delta'}=|\Delta'|$. Since $w^{\Delta}$ and $w^{\widetilde{\Delta}}$ are equal as elements of $\mc{W}$, we have $\ell(w^{\Delta})=\ell(w^{\widetilde{\Delta}})$. But since the subexpressions $w^{\Delta}$ and $w^{\widetilde{\Delta}}$ of $w=s_{i_1}\cdots s_{i_N}$ are reduced, we have $\ell(w^\Delta)=|\Delta|$ and $\ell(w^{\widetilde{\Delta}})=|\widetilde{\Delta}|$; hence $|\Delta|=|\widetilde{\Delta}|$. 

Now we have $|\widetilde{\Delta}|=\height J_\Delta$ and $|\Delta'|=\height J_{\Delta'}$, so that the result now follows from \eqref{almost}.
\end{proof}
\end{proposition}

Yakimov has shown \cite[Theorem 5.7]{Yakimov} that $U_q[w]$ is catenary. We have discussed the uniparameter CGL extension structure of $U_q[w]$. Proposition \ref{crucial inequality for qsc} provides the final condition required for us to apply Theorem \ref{GDME for CGLs with crucial ineq} to $U_q[w]$, giving our main result:
\begin{theorem}
The quantum Schubert cells $U_q[w]$ satisfy the strong Dixmier-Moeglin equivalence.
\end{theorem}

\bibliographystyle{amsplain}
\bibliography{brendansbibliography}\label{references}

\providecommand{\bysame}{\leavevmode\hbox to3em{\hrulefill}\thinspace}
\providecommand{\MR}{\relax\ifhmode\unskip\space\fi MR }
\providecommand{\MRhref}[2]{%
  \href{http://www.ams.org/mathscinet-getitem?mr=#1}{#2}
}
\providecommand{\href}[2]{#2}
\begin{thebibliography}{10}

\bibitem{BRS}
J.~Bell, D.~Rogalski, and S.~J. Sierra, \emph{The {D}ixmier-{M}oeglin
  equivalence for twisted homogeneous coordinate rings}, Israel J. Math.
  \textbf{180} (2010), 461--507.

\bibitem{BCL}
Jason Bell, Karel Casteels, and St{\'e}phane Launois, \emph{Primitive ideals in
  quantum {S}chubert cells: dimension of the strata}, Forum Math. \textbf{26}
  (2014), no.~3, 703--721.

\bibitem{BL}
Jason~P. Bell and St{\'e}phane Launois, \emph{On the dimension of {$H$}-strata
  in quantum algebras}, Algebra Number Theory \textbf{4} (2010), no.~2,
  175--200.

\bibitem{BG}
Ken~A. Brown and Ken~R. Goodearl, \emph{Lectures on algebraic quantum groups},
  Advanced Courses in Mathematics. CRM Barcelona, Birkh\"auser Verlag, Basel,
  2002.

\bibitem{Catoiu}
Stefan Catoiu, \emph{Ideals of the enveloping algebra {$U({\rm sl}_2)$}}, J.
  Algebra \textbf{202} (1998), no.~1, 142--177.

\bibitem{Cauchon}
G{\'e}rard Cauchon, \emph{Effacement des d\'erivations et spectres premiers des
  alg\`ebres quantiques}, J. Algebra \textbf{260} (2003), no.~2, 476--518.

\bibitem{DKP}
C.~De~Concini, V.~G. Kac, and C.~Procesi, \emph{Some quantum analogues of
  solvable {L}ie groups}, Geometry and analysis ({B}ombay, 1992), Tata Inst.
  Fund. Res., Bombay, 1995, pp.~41--65.

\bibitem{Dixmier}
J.~Dixmier, \emph{Id\'eaux primitifs dans les alg\`ebres enveloppantes}, J.
  Algebra \textbf{48} (1977), no.~1, 96--112.

\bibitem{DixmiersBook}
Jacques Dixmier, \emph{Enveloping algebras}, Graduate Studies in Mathematics,
  vol.~11, American Mathematical Society, Providence, RI, 1996, Revised reprint
  of the 1977 translation.

\bibitem{GeigerYakimov}
Joel Geiger and Milen Yakimov, \emph{Quantum {S}chubert cells via
  representation theory and ring theory}, Michigan Math. J. \textbf{63} (2014),
  no.~1, 125--157.

\bibitem{GL}
K.~R. Goodearl and E.~S. Letzter, \emph{Prime and primitive spectra of
  multiparameter quantum affine spaces}, Trends in ring theory ({M}iskolc,
  1996), CMS Conf. Proc., vol.~22, Amer. Math. Soc., Providence, RI, 1998,
  pp.~39--58.

\bibitem{GW}
K.~R. Goodearl and R.~B. Warfield, Jr., \emph{An introduction to noncommutative
  {N}oetherian rings}, second ed., London Mathematical Society Student Texts,
  vol.~61, Cambridge University Press, Cambridge, 2004.

\bibitem{HumphreysCoxeterGroups}
James~E. Humphreys, \emph{Reflection groups and {C}oxeter groups}, Cambridge
  Studies in Advanced Mathematics, vol.~29, Cambridge University Press,
  Cambridge, 1990.

\bibitem{Irving}
Ronald~S. Irving, \emph{Noetherian algebras and nullstellensatz}, S\'eminaire
  d'{A}lg\`ebre {P}aul {D}ubreil 31\`eme ann\'ee ({P}aris, 1977--1978), Lecture
  Notes in Math., vol. 740, Springer, Berlin, 1979, pp.~80--87.

\bibitem{KrauseLenagan}
G{\"u}nter~R. Krause and Thomas~H. Lenagan, \emph{Growth of algebras and
  {G}elfand-{K}irillov dimension}, revised ed., Graduate Studies in
  Mathematics, vol.~22, American Mathematical Society, Providence, RI, 2000.

\bibitem{LLR}
S.~Launois, T.~H. Lenagan, and L.~Rigal, \emph{Quantum unique factorisation
  domains}, J. London Math. Soc. (2) \textbf{74} (2006), no.~2, 321--340.

\bibitem{Lorenz}
Martin Lorenz, \emph{Primitive ideals of group algebras of supersoluble
  groups}, Math. Ann. \textbf{225} (1977), no.~2, 115--122.

\bibitem{LorenzTDeg}
\bysame, \emph{On the transcendence degree of group algebras of nilpotent
  groups}, Glasgow Math. J. \textbf{25} (1984), no.~2, 167--174.

\bibitem{McConnellRobson}
J.~C. McConnell and J.~C. Robson, \emph{Noncommutative {N}oetherian rings},
  revised ed., Graduate Studies in Mathematics, vol.~30, American Mathematical
  Society, Providence, RI, 2001, With the cooperation of L. W. Small.

\bibitem{CauchonMeriaux}
Antoine M{\'e}riaux and G{\'e}rard Cauchon, \emph{Admissible diagrams in
  quantum nilpotent algebras and combinatoric properties of {W}eyl groups},
  Represent. Theory \textbf{14} (2010), 645--687.

\bibitem{Moeglin}
C.~Moeglin, \emph{Id\'eaux primitifs des alg\`ebres enveloppantes}, J. Math.
  Pures Appl. (9) \textbf{59} (1980), no.~3, 265--336.

\bibitem{Yakimov}
Milen Yakimov, \emph{A proof of the {G}oodearl-{L}enagan polynormality
  conjecture}, Int. Math. Res. Not. IMRN (2013), no.~9, 2097--2132.

\bibitem{Zhang}
James~J. Zhang, \emph{On {G}elfand-{K}irillov transcendence degree}, Trans.
  Amer. Math. Soc. \textbf{348} (1996), no.~7, 2867--2899.

\end{thebibliography}

\noindent Jason Bell, University of Waterloo, Department of Pure Mathematics, 200 University Avenue West, Waterloo, Ontario N2L 3G1, Canada\\
E-mail address: {\tt jpbell@uwaterloo.ca}\\

\noindent St\'ephane Launois, University of Kent, School of Mathematics, Statistics, and Actuarial Science, Canterbury, Kent CT2 7NF, UK\\
E-mail address: {\tt S.Launois@kent.ac.uk}\\

\noindent Brendan Nolan, University of Kent, School of Mathematics, Statistics, and Actuarial Science, Canterbury, Kent CT2 7NF, UK\\
E-mail address: {\tt bn62@kent.ac.uk}

\end{document}